\documentclass[11pt,a4paper]{amsart}
\usepackage[latin1]{inputenc}
\usepackage[T1]{fontenc}
\usepackage[english]{babel}

\usepackage{amsmath,amssymb,amsthm,amsfonts}
\theoremstyle{plain}
\newtheorem{thm}{Theorem}

\newtheorem{corollary}[thm]{Corollary}

\theoremstyle{definition}

\newtheorem{rem}[thm]{Remark}
\newtheorem{ex}[thm]{Example}

\newcommand{\extr}{\ensuremath{\mathrm{ext}}}
\newcommand{\cconv}{\overline{\mathrm{co}}}

\newcommand{\C}{\mathbb{C}}
\newcommand{\K}{\mathbb{K}}
\newcommand{\R}{\ensuremath{\mathbb{R}}}
\newcommand{\N}{\ensuremath{\mathbb{N}}}

\newcommand{\tri} {\vert \hspace{-1.5pt} \vert  \hspace{-1.5pt} \vert }

\renewcommand{\leq}{\leqslant}
\renewcommand{\geq}{\geqslant}

\renewcommand{\epsilon}{\varepsilon}

\begin{document}
\begin{center}\small
[Journal of Approximation Theory (2009), \texttt{doi:10.1016/j.jat.2009.07.003}]
\end{center}

\title{On remotality for convex sets in Banach spaces}

\author[Mart\'{\i}n] {Miguel Mart\'{\i}n}
\address[M.~Mart\'{\i}n] {Departamento de An\'{a}lisis Matem\'{a}tico\\
Facultad de Ciencias \\Campus Universitario de Fuentenueva\\
Universidad de Granada \\ E-18071  Granada, Spain,
\textit{E-mail~:} \textit{mmartins@ugr.es}}

\author[Rao]{T.~S.~S.~R.~K.~Rao}
\address[T.~S.~S.~R.~K.~Rao]{Stat--Math Unit\\
Indian Statistical Institute\\
R.~V.~College P.O.\\
Bangalore 560059, India, \textit{E-mail~:}
\textit{tss@isibang.ac.in}}

\thanks{First author partially supported by Spanish MEC project MTM2006-04837 and Junta de
Andaluc\'{\i}a grants FQM-185 and FQM-1438.}

\date{June 18th, 2009}


\maketitle

\begin{abstract} We show that every infinite dimensional Banach space has
a closed and bounded convex set that is not remotal. This settles a
problem raised by Sababheh and Khalil in \cite{S}.
\end{abstract}

%\maketitle

\section{Introduction}

Let $X$ be a real Banach space and let $E \subset X$ be a bounded
set. We write $\extr(E)$ for the set of extreme points of $E$ and
$\cconv(E)$ for the closed (in the norm topology) convex hull of
$E$. If $\tau$ is a locally convex topology in $X$, we will write
$\cconv^{\,\tau}(E)$ to denote the $\tau$-closed convex hull of
$E$. We denote by $B_X$ the closed unit ball of $X$.

The set $E$ is said to be \emph{remotal} from a point $x\in X$, if
there exists a point $e_0 \in E$ such that $D(x,E)=\sup\{\|x-e\|:e
\in E \} = \|x-e_0\|$. The point $e_0$ is  called a \emph{farthest
point} of $E$ from $x$. $E$ is said to be \emph{remotal}
(\emph{densely remotal}) if it is remotal from all (on a dense set)
$x\in X$.  Let $F(x,E)=\{e\in E: D(x,E)=\|x-e\|\}$. In general this
set can be empty. A well known result of Lau (\cite{Lau}) says that
any weakly compact set is densely remotal. It seems to be open, the
question of whether every infinite dimensional Banach space has a
closed and bounded \emph{convex} set that is not remotal. This
question was actually raised in \cite{S} and some partial positive
answers were given in \cite{S} and \cite{R} in the case of reflexive
Banach spaces and Banach spaces that fail the Schur property. The
aim of this note is to give a positive answer to this question. We
follow the notation and terminology of \cite{S} and \cite{R}.

Let us outline the content of this paper. Let $X$ be an infinite
dimensional Banach space and let $X^\ast$ be its topological dual.
Using a classical integral representation theorem, we first show
that $X^\ast$ has a weak$^\ast$-compact convex set $K$ that is not
remotal. This should be compared with \cite[Proposition~1]{DZ} where
the authors exhibited a weak$^\ast$-compact convex set $C \subset
\ell^1$ that has no farthest points. To prove the general result, we
use a stronger form of integral representation theorem for closed
convex bounded sets with the Radon-Nikod\'{y}m property (RNP for
short) due to Edgar (\cite{E}, see \cite[Theorem~16.12]{Phelps}).
Let $E \subset X$ be a weakly closed and bounded set. An interesting
problem that is open is to determine conditions on $\cconv(E)$ so
that  $\cconv(E)$ is remotal from $x$ implies that $E$ is remotal
from $x$. We will give an example showing that $E$ being  norm
closed in a reflexive space is not enough for the validity of
Theorem A in \cite{S}.

\section{Main result}
We first prove a weak$^\ast$-version of \cite[Theorem~A]{S}. In
order to produce a weak$^\ast$-compact convex non-remotal set, it
is enough to show that if $E$ is a weak$^\ast$-compact set having
no vector of maximum length, then the same is true of
$\cconv^{\textrm{weak}^\ast}(E)$ (weak$^\ast$-closed convex hull).
For a compact convex set $K\subset X^\ast$ and for a probability
measure $\mu$ on $K$, let $\gamma(\mu) \in K$ denote its resultant
(or weak integral) with the property
$$
\bigl[\gamma(\mu)\bigr](x)=\int_K k(x)\, d\mu(x) \qquad \bigl(x\in X\bigr).
$$
We refer to \cite{D,Phelps} for the results on integral
representations we use here.

\begin{thm}\label{theorem-1} Let $X$ be an infinite dimensional Banach space.
Let $E \subset X^\ast$ be a weak$^\ast$-closed and bounded set
having no vector of maximum length. Then the weak$^\ast$-closed
convex hull $K$ of $E$ has no vector of maximum length.
Equivalently, if $E$ is not remotal from a point $x\in X$, then
neither is $K$.
\end{thm}

\begin{proof} Let $M= D(0,E)=\sup\{\|e\|:e\in E\}=\sup\{\|k\|:k\in K\}$.
Suppose that there exists $x_0^\ast \in K$ such that
$\|x_0^\ast\|=M$. Let $\mu$ be a probability measure on $K$ with
$\mu(E)=1$ and such that $\gamma(\mu)=x_0^\ast$ (see
\cite[Proposition~1.1]{Phelps}). We fix $\epsilon
>0$ and take $x \in X$  such that $\|x\|=1$ and
$x_0^\ast(x)>M-\epsilon$. Now,
$$
M-\epsilon < x_0^\ast(x)=\int_K x^\ast(x)\, d\mu(x^\ast)=\int_E x^\ast(x)\, d\mu(x^\ast) \leq \int_E
\|x^\ast\|\, d\mu \leq M.
$$
Letting $\epsilon\downarrow 0$, we get that $\displaystyle
\int_E\|x^\ast\|\,d\mu(x^\ast)=M$ and so, $M = \|k\|$ $\mu$-a.e.
Hence $M=\|e\|$ for some $e \in E$. A contradiction. The last part
of the statement is equivalent to the first one just by
translation.
\end{proof}

\begin{corollary}\label{cor-2}
Let $X$ be an infinite dimensional Banach space. Then there exists
a weak$^\ast$-compact convex set $K \subset X^\ast$ that is not
remotal.
\end{corollary}

\begin{proof}
Since $X$ is infinite dimensional, by the well-known
Josefson-Nissenzweig theorem (see \cite[p.~219]{D}), there exists a
sequence $\{x_n^\ast\}_{n \geq 1}$ of unit vectors such that
$x_n^\ast \longrightarrow 0$ in the weak$^\ast$-topology. Consider
the set
$$
E = \left\{\frac{n}{n+1}\,x_n^\ast\ : \ n\in \N\right\}\cup \{0\},
$$
which is clearly a weak$^\ast$-compact set having no vector of
maximum length. Thus, by the above theorem, the weak$^\ast$-closed
convex hull $K$ of $E$ does not have vectors of maximum length, so
$K$ is not remotal from $0$.
\end{proof}

\begin{rem}\label{remark3}
The arguments in Theorem~\ref{theorem-1} and Corollary~\ref{cor-2}
also work in the case of a weakly compact set $E$ and its closed
convex hull $K=\cconv(E)$ (actually, the argument simplifies in
this case and $\epsilon$ is not necessary). Thus, {\slshape in a
Banach space $X$ that fail the Schur property, by taking a
sequence $\{x_n\}_{n \geq 1}$ of unit vectors which converges to
$0$ in the weak topology, we get that the set
$$
K= \cconv\left(\left\{\frac{n}{n+1}\,x_n\ : \ n\in \N\right\}\cup \{0\}\right)
$$
is nonremotal from $0$ (alternatively, the set does not have any
vector of maximal length).}\ This gives an alternative proof of the
main result from \cite{R}.
\end{rem}

\begin{rem} From the above arguments it is easy to see that
{\slshape for a weak$^\ast$-compact set $E \subset X^\ast$ and for
any $x^\ast \in X^\ast$, if the set $F(x^\ast,K)$ of farthest
points in the weak$^\ast$-closed convex hull $K$ of $E$ to
$x^\ast$ is non-empty, then it has a point of $E$.\ } However the
method of proof in \cite{R} has the advantage that it shows that
there is an extreme point of $K$ in $F(x^\ast,K)$. Then by
Milman's theorem \cite[p.~151]{D}, such an extreme point is also
in $E$.
\end{rem}

The following easy example shows that the hypothesis of
weak$^\ast$-closedness can not be omitted on the set $E$ in
Theorem~\ref{theorem-1} (weak-compactness in the case of
Remark~\ref{remark3}).

\begin{ex} {\slshape Let $\{e_n\}_{n\geq 1}$ denote the canonical vector basis
in $\ell^2$. Let $X = \K \oplus_{\infty} \ell^2$, where $\K=\R$ or
$\K=\C$ is the base field and $\oplus_\infty$ means the
$\ell^{\infty}$-direct sum. Consider the set
$$
E=\left\{\left(\frac{n}{n+1},\frac{n}{n+1}\,e_n\right)\ : \ n\in\N\right\}.
$$
Then $E$ is a norm closed set which is not remotal from $0$. Since
$\cconv(E)=\cconv^{\,\textrm{weak}}(E)$ by Mazur's theorem and
$\{e_n\}_{n\geq 1}\longrightarrow 0$ in the weak topology, $(1,0)\in
\cconv(E)$ and so, $\cconv(E)$ is  remotal from $0$.}
\end{ex}

\begin{rem} Let $X$ be a Banach space and let $E \subset X$ be
a weakly closed and bounded set. We do not know if remotality of $K
= \cconv(E)$ from a point always implies that of $E$. Since any
strongly exposed point of $K$ clearly lies in $E$, the answer is
affirmative if the farthest point in $K$ is actually strongly
exposed. We may also ask whether the above question has a positive
answer for RNP sets (see \cite[\S~3]{Bourgin} for these concepts).
\end{rem}

We are now able to present the main result of our paper.

\begin{thm}\label{theorem-main}
Let $X$ be an infinite dimensional Banach space. Then, there exists
a closed and bounded convex set $K$ that is not remotal.
\end{thm}

\begin{proof}
As before, we will construct a closed and bounded set $E$ which is
not remotal from $0$ and show that $K = \cconv(E)$ is also not
remotal from $0$.

In view of Remark~\ref{remark3} (or of \cite{R}), we may assume
without loss of generality that $X$ has the Schur property. Since
$X$ is infinite dimensional, by Rosenthal's $\ell^1$ Theorem (see
\cite[\S~XI]{D}), $X$ contains an isomorphic copy of $\ell^1$. Let
$\tri\cdot\tri$ denote the norm on $X$. Now we will be done if
we can construct in every Banach space $Y=(\ell^1,\tri\cdot\tri)$
isomorphic to $\ell^1$, a closed convex bounded set $K\subseteq Y$
which is not remotal from $0$. Let us write $\tau$ for the
weak$^\ast$-topology of $\ell^1$ as dual of $c_0$ inherited in $Y$.
This is a locally convex topology on $Y$ weaker than the norm
topology and any $\tau$-closed norm-bounded set is compact in this
topology. Observe now that $\tri\cdot\tri$ is not necessarily
weak$^\ast$-lower semi-continuous (i.e.\ $Y$ may not be a dual
space) so, on the one hand, Corollary~\ref{cor-2} does not apply
and, on the other hand, $B_Y$ may not be $\tau$-closed.

Let $\{e_n\}_{n\geq 1}$ be the canonical basis of $\ell^1$.
Consider the set
$$
E =\left\{\frac{n}{n+1}\,\frac{e_n}{\tri e_n\tri }\ : \ n\in \N\right\} \cup \{0\}
\subseteq B_Y
$$
which is $\tau$-compact since $\{e_n\}_{n\geq 1}$ $\tau$-converges
to $0$ and $\tri\cdot\tri$ is equivalent to the usual norm of
$\ell^1$. We consider the set $K = \cconv^{\,\tau}(E) \subset Y$,
which is $\tau$-compact since it is $\tau$-closed and norm-bounded
(indeed, $E$ is contained in the $\tau$-closed set $M B_{\ell^1}$
for some $M>0$, so $K\subset M B_{\ell^1}$).

\noindent{\emph{Claim.}}\ $K\subseteq B_Y$.\ Indeed, since
$\ell^1$ (and so $Y$) has the RNP, $K$ is a set with the RNP.
Therefore, we have $K = \cconv\bigl(\extr(K)\bigr)$ (closure in
norm, see \cite[\S~3]{Bourgin}). As $K$ and $E$ are
$\tau$-compact, Milman's theorem gives us that $\extr(K) \subseteq
E$ (see \cite[p.~151]{D}). Therefore, we have
$$
K=\cconv\bigl(\extr(K)\bigr) \subseteq \cconv(E)\subseteq
\cconv^{\,\tau}(E)=K,
$$
so $K=\cconv(E)\subseteq B_Y$ as claimed.

Suppose $K$ is remotal from $0$ in $Y$. As $D(0,E)=1$ and
$K\subseteq B_Y$, we also have $D(0,K)=1$. Therefore, there is a
vector $y_0\in K$ with $\tri y_0\tri=1$, and we may pick a
functional $y_0^\ast\in Y^\ast$ with
$$
\tri y_0^\ast\tri=1 \qquad \text{and} \qquad y_0^\ast(y_0)=1.
$$
As $K$ is a separable closed convex bounded set with the RNP,
Edgar's integral representation theorem (\cite{E}, see
\cite[Theorem~16.12]{Phelps}), gives us that there exists a
probability measure $\mu$ on $K$ with $\mu\bigl(\extr(K)\bigr)=1$
(so $\mu(E)=1$) such that
$$
1=y_0^\ast(y_0)=\int_{K} y_0^\ast(y)
\,d\mu(y)= \int_{E} y_0^\ast(y)
\,d\mu(y) \leq \int_{E}\tri y\tri
\,d\mu(y)\leq 1.
$$
Therefore, $\tri y \tri=1$ $\mu$-a.e.\ in $E$, which is clearly
false. Thus we get a contradiction and $K$ is nonremotal from $0$.
\end{proof}

Since remotality from $0$ is equivalent to having a vector of
maximal norm, we get the following corollary.

\begin{corollary}
Let $X$ be an infinite-dimensional Banach space. Then there is a
closed convex set $K$ contained in the \emph{open} unit ball of
$X$ such that $\sup\{\|x\|\,:\,x\in K\}=1$.

\end{corollary}

\begin{rem} Similar to Remark~4 (see also Remark~6), let us note that
{\slshape for a separable weakly closed and bounded set $E$ such
that its closed convex hull $K$ has the RNP, our arguments show
that if $F(x,K) \neq \emptyset$ then it has an extreme point of
$K$.}
\end{rem}

\begin{rem}
{\slshape Going into the details of the proofs of
Remark~\ref{remark3} and Theorem~\ref{theorem-main}, one realizes
that for every infinite-dimensional Banach space $X$, there is a
locally convex Hausdorff topology $\tau$, which is weaker than the
norm topology and  such that there is a $\tau$-compact convex set
$K$ which is not remotal (from $0$).\ } Indeed, if $X$ does not have
the Schur property, then the set $K$ is actually weak compact.
Otherwise, $X$ contains a subspace $Y$ isomorphic to $\ell^1$, and
the set $K\subset Y$ is compact for the topology $\tau'$ of $Y$
which it inherits from the weak$^\ast$ topology of $\ell^1$ as dual
of $c_0$. Since we may extend the topology $\tau'$ of $Y$ to a
locally convex Hausdorff topology $\tau$ of $X$ (still weaker than
the norm topology of $X$), we get that $K$ is $\tau$-compact, as
desired.
\end{rem}

\section*{Acknowledgement} This research was done during
a visit of first author to the Bangalore centre of the Indian
Statistical Institute, where he participated in the \emph{Workshop
on Geometry of Banach spaces} funded by the National Board for
Higher Mathematics (NBHM), Government of India. He also would like
to thanks the people at the ISI Bangalore for the warm hospitality
received there.


\begin{thebibliography}{99}

\bibitem{Bourgin} \textsc{R.~R.~Bourgin}, \emph{Geometric
    Aspects of Convex Sets with the Radon-Nikodym Property}, Lecture
    Notes in Math. \textbf{993}, Springer-Verlag, Berlin 1983.

  \bibitem{DZ} \textsc{R.~Deville and V.~E.~Zizler}, {\it Farthest points in
  $w\sp \ast$-compact sets},
   Bull. Austral. Math. Soc. \textbf{38} (1988), 433--439

  \bibitem{D} \textsc{J.~Diestel}, {\em Sequences and series in Banach spaces},
   Graduate Texts in Mathematics, {\bf 92}. Springer-Verlag, New York, 1984.

\bibitem {E} \textsc{G.~A.~Edgar}, {\it Extremal integral
    representations},
    J. Funct. Anal. \textbf{23} (1976), 145--161.

\bibitem{Lau} \textsc{K.-S.~Lau}, {\em Farthest Points in Weakly
    Compact Sets}, Israel J. Math. {\bf 22} (1975), 168--174.

\bibitem{Phelps} \textsc{R.~R.~Phelps}, {\em Lectures on
    Choquet's Theorem}, Second Edition, Lecture Notes in Math.
    {\bf 1757}, Springer, Berlin, 2001.

\bibitem{R} \textsc{T.~S.~S.~R.~K.~Rao}, {\em Remark on a paper of
    Sababheh and Khalil},  Numerical Functional Analysis and Optimization \textbf{30} (2009), 822--824.

\bibitem{S} \textsc{M.~Sababheh and R.~Khalil}, {\em Remotality of
    closed bounded convex sets}, Numerical Functional Analysis and
    Optimization {\bf 29} (2008), 1166-1170.

\end{thebibliography}
\end{document}